\documentclass[12pt,a4paper,twoside]{article}
\usepackage[english]{babel}
\usepackage{amssymb}
\usepackage{inputenc}
\usepackage{graphicx}
\newtheorem{thm}{Theorem}[section]
\newtheorem{lem}{Lemma}[section]

\newtheorem{Def}{Definition}[section]

\newtheorem{Rem}{Remark}[section]

\def\bc{\begin{center}}
\def\ec{\end{center}}
\date{}
\begin{document}

\author{S. Albeverio $^{1},$ B.A.
Omirov  $^2,$ A.Kh. Khudoyberdiyev$^3$}

\title{\bf On the classification of complex Leibniz superalgebras with characteristic sequence
$(n-1, 1 | m_1, \dots, m_k)$ and nilindex $n+m$. }

\maketitle

\begin{abstract} In this work we investigate the complex Leibniz superalgebras with
characteristic sequence $(n-1, 1 | m_1, \dots, m_k)$ and with
nilindex equal to $n+m.$ We prove that such superalgebras with the
condition $m_2\neq0$ have nilindex less than $n+m$. Therefore the
complete classification of Leibniz algebras with characteristic
sequence $(n-1, 1 | m_1, \dots, m_k)$ and with nilindex equal to
$n+m$ is reduced to the classification of filiform Leibniz
superalgebras of nilindex equal to $n+m,$ which was provided in
\cite{AOKh} and \cite{GKh}.
\end{abstract}

\medskip
$^1$ Institut f\"{u}r Angewandte Mathematik, Universit\"{a}t Bonn,
Wegelerstr. 6, D-53115 Bonn (Germany); SFB 611, BiBoS; CERFIM
(Locarno); Acc. Arch. (USI),  e-mail: \emph{albeverio@uni-bonn.de}

$^2$ Institute of Mathematics and Information Technologies,
Uzbekistan Academy of Science, F. Hodjaev str. 29, 100125,
Tashkent (Uzbekistan), e-mail: \emph{omirovb@mail.ru}

 $^{3}$ Institute of Mathematics and Information Technologies,
Uzbekistan Academy of Science, F. Hodjaev str. 29, 100125,
Tashkent (Uzbekistan), e-mail: \emph{khabror@mail.ru}



\section{Introduction}

During many years the theory of Lie superalgebras has been
actively studied by many mathematicians and physicists. A
systematic exposition of basic Lie superalgebras theory can be
found in \cite{Kac}. Many works have been devoted to the study of
this topic, but unfortunately most of them do not deal with
nilpotent Lie superalgebras. In works \cite{Rosa}, \cite{Gilg},
\cite{GKh} the problem of the description of some classes of
nilpotent Lie superalgebras has studied. It is well known that Lie
superalgebras are a generalization of Lie algebras. In the same
way, the notion of Leibniz algebras, which were introduced in
\cite{L}, can be generalized to Leibniz superalgebras. The
elementary properties of Leibniz superalgebras were obtained in
\cite{AAO2}. For nilpotent Leibniz superalgebras the description
of the case of maximal nilindex (nilpotent Leibniz superalgebras
distinguished by the feature of being singly-onegenerated) is not
difficult and was done in \cite{AAO2}. However, the next stage
(the description of Leibniz superalgebras with dimensions of the
even and odd parts equal to n and m, respectively, and of nilindex
$n+m$) is a very problematic one. It should be noted that such Lie
superalgebras were classified in \cite{GKh}. Due to the great
difficulty of solving in general the problem of description of
Leibniz superalgebras of nilindex $n+m$, some restrictions on the
characteristic sequence should be added, in particular, since the
graded anticommutative identity does not hold in non-Lie Leibniz
superalgebras. In the description of the structure of Leibniz
superalgebras the crucial task is to prove the existence of a
suitable basis (the so-called adapted basis) in which the table of
multiplication of the superalgebra has the most convenient form.
In the present paper we investigate the Leibniz superalgebras with
the characteristic sequence $C(L)=(n-1,1 | m_1, m_2, \dots, m_k)$
and with nilindex equal to $n+m.$ Actually, the classification of
such superalgebras in the case where $C(L)=(n-1, 1 | m)$ was
obtained in \cite{AOKh} and the main result of the present work
consist of the following fact: Leibniz superalgebras with
characteristic sequences equal to $C(L)=(n-1,1 | m_1, m_2, \dots,
m_k)$ ($m_2\neq0$) have nilindex less than $n+m.$ Therefore, the
classification of Leibniz superalgebras of nilindex $n+m$ and with
characteristic sequence equal to $(n-1,1 | m_1, m_2, \dots, m_k)$
is reduced to the case where the
 characteristic sequence is equal to $(n-1, 1 | m)$, and this case, as we
mentioned above, has already been solved in \cite{AOKh}. In this
way we made further step in the solution of the problem of the
classification of complex Leibniz superalgebras of nilindex $n+m.$

Throughout this work we shall consider spaces and (super)algebras
over the field of complex numbers.

\section{Preliminaries}
We recall the notions of Lie and Leibniz superalgebras.

\begin{Def} \emph{\cite{Kac}} {\em A $Z_2$-graded vector space
$G=G_0\oplus G_1$ is called} a Lie superalgebra {\em if it is
equipped with a product $[-,-]$ which satisfies the following
conditions:

1. $[G_\alpha,G_\beta]\subseteq G_{\alpha+\beta(mod\ 2)}$  for any
$\alpha,\beta\in Z_2,$

2. $[x,y]=-(-1)^{\alpha\beta}[y,x],$ for any $x\in G_\alpha,$
$y\in G_\beta,$

3. $(-1)^{\alpha\gamma}
[x,[y,z]]+(-1)^{\alpha\beta}[y,[z,x]]+(-1)^{\beta\gamma}[z,[x,y]]=0$}
---Jacobi superidentity, \\ {\em for any $x\in G_\alpha,$ $y\in G_\beta,$
$z\in G_\gamma, \ \alpha, \beta, \gamma\in\mathbb{Z}_2.$}
\end{Def}

\begin{Def} \emph{\cite{AAO2}} {\em A $Z_2$-graded vector space $L=L_0\oplus
L_1$ is called} a Leibniz superalgebra {\em if it is equipped with
a product $[-, -]$ which satisfies the following conditions:

1. $[L_\alpha,L_\beta]\subseteq L_{\alpha+\beta(mod\ 2)}$  for any
$\alpha,\beta\in Z_2,$

2.  $[x, [y, z]]=[[x, y], z] - (-1)^{\alpha\beta} [[x, z], y]$}
---Leibniz superidentity, {\em for any $x\in L,$ $y \in L_\alpha,$ $z \in
L_\beta$.}
\end{Def}
Evidently, the subspaces $G_0$ and $L_0$ are Lie and Leibniz
algebras, respectively.

It should be noted that if in a Leibniz superalgebra L the
identity:
  $$
[x,y]=-(-1)^{\alpha\beta}[y,x],
$$
holds for any $x\in L_\alpha$ and $y\in L_\beta,$ then the Leibniz
superidentity can easily be transformed into the Jacobi
superidentity. Thus, Leibniz superalgebras are generalizations of
both Lie (super)algebras and Leibniz algebras. For examples of
Leibniz superalgebras we refer to \cite{AAO2}.

The set of Leibniz superalgebras with dimensions of the even part
$L_0$ and the odd part $L_1$, respectively equal to $n$ and $m,$
shall be denoted by $Leib_{n,m}.$

For a given Leibniz superalgebra $L$ we define a descending
central sequence in the following way:
$$
L^1=L,\quad L^{k+1}=[L^k,L^1], \quad k \geq 1.
$$

\begin{Def} {\em A Leibniz superalgebra $L$ is called}
nilpotent, {\em if there exists  $s\in\mathbb N$ such that
$L^s=0.$ The minimal number $s$ with this property is called}
index of nilpotency (or nilindex) of the superalgebra $L.$
\end{Def}

\begin{Def} {\em The set ${\mathcal R}(L)=\left\{ x\in L\ |\ [y,
x]=0 \ \mbox{for any} \  y\in L \right\}$ is called} the right
annihilator of a superalgebra $L$.
\end{Def}

Using the Leibniz superidentity it is not difficult to see that
${\mathcal R}(L)$ is an ideal of the superalgebra $L$. Moreover,
elements of the form $[a,b]+(-1)^{ab}[b,a]$ belong to ${\mathcal
R}(L)$.

The following theorem describes nilpotent Leibniz superalgebras
with maximal nilindex.

\begin{thm} \label{t1} \emph{\cite{AAO2}} Let $L=L_0\oplus L_1$  be a Leibniz superalgebra
from $Leib_{n,m }$ with nilindex equal to $n+m+1.$ Then $L$ is
isomorphic to one of the following two non isomorphic
superalgebras:
$$
[e_i,e_1]=e_{i+1},\ 1\le i\le n-1;\quad \left\{ \begin{array}{ll} [e_i,e_1]=e_{i+1},& 1\le i\le n+m-1, \\
{[}e_i,e_2{]}=2e_{i+2}, & 1\le i\le n+m-2.\\ \end{array}\right.
$$
where omitted products are equal to zero and $\{e_1, e_2, \dots,
e_n\}$ is the basis of the superalgebra $L$.
\end{thm}

\begin{Rem} \label{Rem1} {\em From the description of Theorem \ref{t1} we have that if the odd part
$L_1$ of the superalgebra $L$ is non trivial, then either $m=n$ or
$m=n+1$ and the table of multiplication of the second superalgebra
in a graded basis $\{x_1, \dots, x_n, y_1, \dots, y_m\}$ can be
written in the following form:}
$$[y_j, x_1] = y_{j+1}, \ 1 \leq j \leq m-1, \quad [x_1, y_1] = \frac 1 2 y_2,$$
$$[x_i, y_1] = \frac 1 2 y_i, 2 \leq i \leq m, \quad \quad \ \ \ [y_1, y_1] = x_1,$$
$$[y_j, y_1] = x_{j+1}, \quad 2 \leq j \leq n-1,$$
$$[x_i,x_1]=x_{i+1}, \quad 1 \leq i \leq n-1.$$
\end{Rem}

Let $L=L_0\oplus L_1$ be a nilpotent Leibniz superalgebra. For an
arbitrary element $x\in L_0,$ the operator of right multiplication
$R_x$ (defined as $R_x :y \rightarrow [y,x]$) is a nilpotent
endomorphism of the space $L_i,$ where $i\in \{0, 1\}.$ Denote by
$C_i(x)$ ($i\in \{0, 1\}$) the descending sequence of the
dimensions of Jordan blocks of the operator $R_x.$ Consider the
lexicographical order on the set $C_i(L_0)$.

\begin{Def} {\em The sequence
$$
C(L)=\left( \left.\max\limits_{x\in L_0\setminus [L_0,L_0]}
C_0(x)\ \right|\ \max\limits_{\widetilde x\in L_0\setminus
[L_0,L_0]} C_1\left(\widetilde x\right) \right)
$$
is said to be} the characteristic sequence of the Leibniz
superalgebra $L.$
\end{Def}

Similarly to \cite{Gilg} (corollary 3.0.1) it can be proved that
the characteristic sequence is invariant under isomorphisms.

Further we need the following definition.
\begin{Def} {\em A Leibniz algebra $L$ of dimension $n$ is said to be} filiform {\em if
$dim L^i=n-i$ for $2\leq i \leq n.$}
\end{Def}

\begin{lem}\label{lem1} \emph{\cite{AO}} Let $L$ be an $n$-dimensional Leibniz algebra. Then the
following statements are equivalent: $$ \quad \quad  a) \ C(L) =
(n - 1, 1);$$
$$\quad \quad  b)  \ L \ \mbox{is a filiform Leibniz algebra};$$
$$ c) \  L^{n-1}\neq0 \ \mbox{and} \  L^n =0.$$
\end{lem}

Let $L$ be a Leibniz superalgebra from $Leib_{n,m}$ with
characteristic sequence equal to $(n-1,1 | m_1, m_2, \dots, m_k),$
(where $m_1+m_2+\dots+m_k=m).$ Since in \cite{GKh} the Leibniz
superalgebras with characteristic sequence and nilindex equal to
$(n-1, n |m)$ and $n+m,$ respectively, have already been obtained,
we shall henceforth reduce our investigation to the case where
$m_2\neq0.$

From Lemma \ref{lem1} we can conclude that the even part $L_0$ of
$L$ is a filiform Leibniz algebra. Due to the description of
filiform Leibniz algebras in \cite{GJKh}, \cite{GO}, \cite{Ve} we
can obtain the existence of an adapted basis in superalgebra with
$C(L)=(n-1, 1 | m_1, m_2, \dots, m_k),$ according to the following
theorem

\begin{thm}\label{t2} Let $L=L_0\oplus L_1$ be a superalgebra from $Leib_{n,m}$
with characteristic sequence equal to $(n-1, 1 | m_1, m_2, \dots,
m_k).$ Then there exists a basis $\{ x_1, x_2, \ldots, x_n, y_1,
y_2, \ldots, y_m\}$ of $L,$ in which the multiplication satisfies
one of the following three conditions:

a) $[x_1,x_1]=x_3,$

\ \quad $[x_i,x_1]=x_{i+1},$ $2\le i\le n-1,$

\ \quad $[y_j,x]=y_{j+1},$ $j\notin \{m_1, m_1+m_2, \dots,
m_1+m_2+\dots+m_k\},$

\ \quad $[y_j,x]=0, \ \quad j\in \{m_1, m_1+m_2, \dots,
m_1+m_2+\dots+m_k\},$ \\for some $x\in L_0\setminus L_0^2,$

 \ \quad
$[x_1,x_2]=\alpha_4x_4+\alpha_5x_5+\ldots+\alpha_{n-1}x_{n-1}+\theta
x_n,$

\ \quad $[x_j,x_2]=\alpha_4x_{j+2}+\alpha_5x_{j+3}+\ldots+\alpha_{n+2-j}x_{n},$ $2\le j\le n-2,$\\
where the omitted products in $L_0$ are equal to zero;

b) $[x_1,x_1]=x_3,$

\ \quad $[x_i,x_1]=x_{i+1},$ $3\le i\le n-1,$

\ \quad $[y_j,x]=y_{j+1},$ $j\notin \{m_1, m_1+m_2, \dots,
m_1+m_2+\dots+m_k\},$

\ \quad $[y_j,x]=0, \ \quad j\in \{m_1, m_1+m_2, \dots,
m_1+m_2+\dots+m_k\},$ \\for some $x\in L_0\setminus L_0^2,$

\ \quad $[x_1,x_2]=\beta_4x_4+\beta_5x_5+\ldots+\beta_{n}x_{n},$

\ \quad $[x_2,x_2]=\gamma x_{n},$

\ \quad $[x_j,x_2]=\beta_4x_{j+2}+\beta_5x_{j+3}+\ldots+\beta_{n+2-j}x_{n},$ $3\le j\le n-2,$\\
where the omitted products in $L_0$ are equal to zero;

c) $[x_i,x_1]=x_{i+1},$ $2\le i\le n-1,$

\ \quad $[x_1,x_i]= - x_{i+1},$ $3\le i\le n-1,$

\ \quad $[y_j,x]=y_{j+1},$ $j\notin \{m_1, m_1+m_2, \dots,
m_1+m_2+\dots+m_k\},$

\ \quad $[y_j,x]=0, \ \quad j\in \{m_1, m_1+m_2, \dots,
m_1+m_2+\dots+m_k\},$ \\for some $x\in L_0\setminus L_0^2,$

\ \quad $[x_1,x_1]=\theta_1x_n,$

\ \quad $[x_1,x_2]=-x_3+\theta_2x_n,$

\ \quad $[x_2,x_2]=\theta_3x_n,$

 \ \quad
$[x_i,x_j]=-[x_j,x_i]\in lin\langle x_{i+j+1}, x_{i+j+2}, \ldots,
x_n\rangle,$ $2\le i < j\le n-2.$
\end{thm}

\section{On the classification of Leibniz superalgebras with characteristic sequence $(n-1, 1 | m_1, m_2, \dots, m_k)$
and nilindex $n+m$ ($m_2\neq 0$).}

Let $L$ satisfy to the conditions of Theorem \ref{t2} and let
$\{x_1, x_2, \dots, x_n, y_1, y_2, \dots, y_m\}$ be the adapted
basis of $L.$ It is not difficult to see that if $L$ has nilindex
equal to $n+m,$ then the superalgebra $L$ has two generators (due
to Theorem \ref{t1} we have a description of singly-generated
Leibniz superalgebras, which have nilindex $n+m+1$) and
$dimL^{i}=n+m-2$ for $2\leq i \leq n+m.$ It should be noted that
the filiform Leibniz algebra $L_0$ has also two generators, $x_1$
and $x_2$.

\begin{lem}\label{l1} In three the classes of superalgebras of Theorem \ref{t2}
instead of the element $x$ one can choose the element $x_1.$
\end{lem}

\textit{Proof.} Without loss of generality we can assume that $x$
has the form: $x=A_1x_1+A_2x_2,$ where $(A_1,A_2)\ne (0,0).$

Let us consider the first class of superalgebras of Theorem
\ref{t2} and investigate three cases.

{\bf Case 1.} Let $A_1(A_1+A_2) \ne 0.$ Then applying the
following change of basis:
$$
x_1'=A_1x_1+A_2x_2,\quad
x_2'=(A_1+A_2)x_2+A_2(\theta-\alpha_n)x_{n-1},
$$ $$
x_i'=[x_{i-1}',x_1'], \quad 3\le i\le n, \quad y_j'=y_j,\quad 1\le
j\le m,
$$
we obtain that the first four multiplications in class a) do not
change.

{\bf Case 2.} Let $A_1=0.$ Let us make a change of basis as
follows:

$$ x_1'=x_1+aA_2x_2,\mbox{ where } a(1+aA_2)\ne
0,$$
$$ \ x_2'=(1+A_2)x_2+aA_2(\theta- \alpha_n) x_{n-1}, \quad  x_i'=[x_{i-1}',x_1'], \ 3\le i\le
n,$$
$$y_j'=y_j, \qquad  \quad j\in \{1, m_1+1, m_1+m_2+1, \dots,
m_1+m_2+\dots+m_{k-1} +1\},$$
 $$y_j'=[y_{j-1}',x_1'], \ j\notin \{1, m_1+1, m_1+m_2+1, \dots, m_1 + m_2 + \dots +
m_{k-1} +1\}.$$

If we choose a sufficiently big value of the parameter $a$ then we
obtain that the first four multiplications in the class a) also do
not change. Indeed, the first three multiplications do not change
by the construction and the products $[y_j',x_1']$ for $j\in
\{m_1, m_1+m_2, \dots, m_1+m_2+\dots+m_k\}$ are equal to zero,
because otherwise we easily can get a contradiction with the
characteristic sequence or nilpotence conditions.

{\bf Case 3.} Let $A_1\ne 0$ and $A_1=-A_2.$ Then taking the
following transformation of basis:
$$
x_1'=A_1x_1-A_1x_2+ax_2,\quad x_2'=ax_2+(a-A_1)(\theta -\alpha_n
)x_{n-1}, \ (a\ne 0),
$$ $$
x_i'=[x_{i-1}', x_1'],\ 3 \le i \le n, \ y_j'=y_j, \ j\in \{1,
m_1+1, m_1+m_2+1, \dots, m_1+m_2+\dots+m_{k-1} +1\},$$
 $$y_j'=[y_{j-1}',x_1'], \ j\notin \{1, m_1+1, m_1+m_2+1, \dots, m_1 + m_2 + \dots +
m_{k-1} +1\},$$ it is not difficult to check that for sufficiently
small values of the parameter $a$ the first four multiplications
in the class a) are preserved.

Thus, we have shown that in the first case of superalgebras of
Theorem \ref{t2} instead of element $x$ one can choose element
$x_1.$

Let us consider the class b) of Theorem \ref{t2}.

If $A_1\neq0,$ then applying a transformation of basis of the
form:
$$x_1'=A_1x_1+A_2x_2, \ x_2'=x_2-\frac{A_2 \gamma}{A_1}x_{n-1},$$ $$x_3'=[x_1',x_1'], \ x_i'=[x_{i-1}',x_1'],
\ 4\le i\le n, \ y_j'=y_j, \ 1\le j\le m,
$$ we obtain that the first four multiplications do not change.

If $A_1=0,$ then the following change of basis:
$$x_1'=x_1+aA_2x_2, \ (a\neq0), \ x_2'=x_2-aA_2\gamma x_{n-1},$$ $$x_3'=[x_1',x_1'], \ x_i'=[x_{i-1}',x_1'],
\ 4\le i\le n,$$
$$y_j'=y_j, \qquad \quad j\in \{1,
m_1+1, m_1+m_2+1, \dots, m_1+m_2+\dots+m_{k-1} +1\},$$
 $$y_j'=[y_{j-1}',x_1'], \ j\notin \{1, m_1+1, m_1+m_2+1, \dots, m_1 + m_2 + \dots +
m_{k-1} +1\},$$ with a sufficiently big value of the parameter $a$
allow to conclude that first four multiplications in class b) do
not change.

Now consider the class c) of Theorem \ref{t2}.

If $A_1\neq0,$ then applying a transformation of basis of the
form:
$$x_1'=A_1x_1+A_2x_2, \ x_2'=x_2, \ x_i'=[x_{i-1}',x_1'],
\ 3\le i\le n, \ y_j'=y_j, \ 1\le j\le m,$$ we obtain that the
first four multiplications are preserved.

If $A_1=0,$ then take the transformation of basis:

$$x_1'=x_1+aA_2x_2, \ (a\neq0), \ x_2'=x_2, \ x_i'=[x_{i-1}',x_1'],
\ 3\le i\le n,$$
$$y_j'=y_j, \qquad \quad j\in \{1,
m_1+1, m_1+m_2+1, \dots, m_1+m_2+\dots+m_{k-1} +1\},$$
 $$y_j'=[y_{j-1}',x_1'], \ j\notin \{1, m_1+1, m_1+m_2+1, \dots, m_1 + m_2 + \dots +
m_{k-1} +1\}.$$ Then choosing a sufficiently big value of the
parameter $a$ allow us to conclude that the first four products in
the case c) of Theorem \ref{t2} do not change.

Thus, we have proven that in the three classes of superalgebras of
Theorem \ref{t2} instead of the element $x$ we can choose element
$x_1$. \hfill $\square$

Since the superalgebra $L=L_0\oplus L_1$ has two generators the
possible cases are as follow: both generators lie in $L_0;$ one
generator lies in $L_0$ and the another one lies in $L_1;$ both
generators lie in $L_1.$

We shall not consider the case where both generators lie in the
even part (since $m_2\neq0$). Firstly we consider the second
possible, i.e. the case where one of the generators lies in $L_0$
and the another one lies in $L_1.$ It is easy to see that there
exist some $m_j, \ 0 \leq j \leq k-1$ (here $m_0 = 0$), such that
$y_{m_1 + m_2 + \dots + m_j+1} \notin L^2.$ By a shifting of basic
elements one can assume that $m_j=m_0$, i.e. the basic element
$y_1\notin L^2.$ Of course, by this choice the condition $m_1 \geq
m_2 \dots m_k,$ is broken, but we shall not use this condition in
our study further. Thus, as generators we can choose the elements
$A_1x_1+A_2x_2$ and $y_1.$

Let us introduce the notations
$$[x_i, y_1] = \sum\limits_{j=2}^m\alpha_{i,j}y_j, \ 1 \leq i \leq n,
\ [y_s, y_1] = \sum\limits_{t=1}^n\beta_{s,t}x_t, \ [y_p, x_2]
=\sum\limits_{q=2}^{m}\gamma_{p,q}y_q, \ 1 \leq s, p \leq m.$$

\begin{thm}\label{t3} Let $L$ be a Leibniz superalgebra from $Leib_{n,m}$ with characteristic sequence
$Ñ(n-1,1 | m_1, m_2, \dots, m_k),$ where $m_1\geq2, \ n\geq 4$.
Let the elements $A_1x_1+A_2x_2, \ y_1$ be generators and $x_1 \in
L^2.$ Then $L$ has a nilindex less than $n+m.$
\end{thm}
\textit{Proof.} Since $x_1 \in  L^2,$ then $x_2\notin L^2$ and
therefore as a generator of the $L$ which lies in $L_0$ we can
choose $x_2.$ Let us assume the contrary, i.e. the nilindex of the
superalgebra $L$ is equal to $n+m.$ Then we have
$$ L = \{x_1, x_2, \dots, x_n, y_1, y_2, \dots, y_m\}, \ \  L^2 = \{x_1, x_3, \dots, x_n, y_2, y_3, \dots, y_m\}.$$
Since $x_1$ is generator in the filiform Leibniz algebra $L_0$
then it should lie in linear span of products $[y_s,y_1], \ 1 \leq
s \leq m.$ Therefore, there exists some $s$ ($1 \leq s \leq m$)
such that $\beta_{s,1}\neq 0.$

Denote by $y_{m_1+m_2+ \dots +m_{t_0}+1}$ the basic element which
is the earliest generated among the elements $\{y_{m_1+1},
y_{m_1+m_2+1}, \dots, y_{m_1+m_2+ \dots +m_{k-1}+1}\},$ i.e.
$y_{m_1+m_2+ \dots +m_{t_0}+1}$ is the basic element which first
is absent among the elements $\{y_{m_1+1}, y_{m_1+m_2+1}, \dots,
y_{m_1+m_2+ \dots +m_{k-1}+1}\}$ in the descending central
sequence. Therefore, the element $y_{m_1+m_2+ \dots +m_{t_0}+1}$
is generated by products of elements of the form either $[x_i,
y_1], 2 \leq i \leq n$ or $[y_j, x_2], 1 \leq j \leq m_1.$

Let us show that $x_1 \notin L^3.$ Indeed, if $x_1$ lies in $L^3,$
then it should be generated by the product $[y_{m_1+m_2+ \dots
m_{t_0}+1}, y_1]$. From the nilpotence condition we have that
$x_1$ is generated by the products $[[x_2, y_1], y_1],$ $[[y_1,
x_2], y_1]$ (since the basic elements $x_3, x_4, \dots, x_n$ and
$y_2, y_3, \dots, y_{m_1}$ are obtained by the products involving
 $x_1$). Hence for generating the element $y_{m_1+m_2+ \dots
+m_{t_0}+1}$ it is enough to consider the cases where $i = 2, \ j
= 1.$

From the equalities
$$[[x_2, y_1], y_1]=\frac{1}{2}[x_2, [y_1, y_1]]=\frac{1}{2}[x_2, \sum\limits_{s=1, s\neq2}^n\beta_{1,s}x_s] =
\sum\limits_{t\geq 3}(*)x_t$$ (where by the symbol $(*)$ we denote
the coefficients of the basic elements $x_t$), we have that the
element $x_1$ is not present in the decomposition of the product
$[[x_2, y_1], y_1]$. If $x_1$ is generated from the product
$[[y_1, x_2], y_1],$ then the expression $[[y_1, x_2], y_1] +
[[x_2, y_1], y_1]$ lies in $ {\mathcal R}(L)$ and $x_1$ appear in
its decomposition. Using the table of multiplication in the
algebra $L_0$ from Theorem \ref{t2} we establish that multiplying
the expression $[[y_1, x_2], y_1] + [[x_2, y_1], y_1]$ on the
right side by the element $x_1$ sufficiently many times we obtain
$x_n\in{\mathcal R}(L).$ Then repeating this procedure finally we
obtain $x_1\in{\mathcal R}(L).$ Thus, we obtain a contradiction,
because $x_1 \notin {\mathcal R}(L).$ Therefore, $x_1 \notin L^3.$

Let us consider the cases a) and b). The condition $x_1 \notin
L^3$ leads to $\beta_{1,1} \neq 0.$ Since in the cases a) and b)
the elements $x_i, \ 3\leq i \leq n$ lie in ${\mathcal R}(L),$
then we can put $x'_1 = \beta_{1,1}x_1 + \beta_{1,3}x_3 + \dots +
\beta_{1,n}x_n$ and suppose that $[y_1, y_1] = x_1.$

Consider the subsuperalgebra generated by $y_1$ ($<y_1>$). Then it
is easy to check that $<y_1> = \{x_1, x_3, \dots, x_n, y_1, y_2,
\dots, y_{m_1}\}.$ Since this subsuperalgebra is single-generated
then from Remark \ref{Rem1} we have
$$[y_j, x_1] = y_{j+1}, \ 1 \leq j \leq m_1-1, \ \ [x_1, y_1] = \frac 1 2 y_2,$$
$$[x_i, y_1] = \frac 1 2 y_i, 3 \leq i \leq m_1, \ \ [y_1, y_1] = x_1,$$
$$[y_j, y_1] = x_{j+1}, \quad 2 \leq j \leq n-1,$$
where if $n+m_1$ is even then $m_1=n-1$ and if $n+m_1$ is odd then
$m_1=n.$

From above products we have $[y_1, x_1] + [x_1, y_1] =
\frac{3}{2}y_2$ which yields $y_i\in {\mathcal R}(L)$ for $2\leq i
\leq m_1.$

Using the fact that in the cases a) and b) the product $[x_1,x_2]$
belongs to ${\mathcal R}(L)$ and the following equalities hold
$$ [y_i, x_2] = [[y_{i-1}, x_1], x_2] = [y_{i-1}, [x_1, x_2]] + [[y_{i-1},
x_2], x_1]= [[y_{i-1}, x_2], x_1]=$$ $$\gamma_{1,2}y_{i+1} +
(*)y_{i+2}+\dots + (*)y_{m_1}+(*)y_{m_1+2} + \dots +
(*)y_{m_1+m_2}+ (*)y_{m_1+m_2+2}+ \dots +(*)y_m,$$ for $2\leq i
\leq m_1,$ we can conclude that $y_{m_1+m_2+ \dots +m_{t_0}+1}
\notin L^3,$ i.e. $dimL^3< n+m-3.$ Hence we obtain a contradiction
with the assumption that the nilindex is equal to $n+m.$

Now, consider superalgebras from the class b) of Theorem \ref{t2}.
We have $$L^2 = \{x_1, x_3, \dots, x_n, y_2, y_3, \dots, y_m\},$$
$$L^3 = \{x_3, \dots, x_n, y_2, y_3, \dots, y_m\}.$$

Consider the equalities
$$[y_1, x_3] = [y_1, [x_2, x_1]] = [[y_1,
x_2], x_1] - [[y_1, x_1], x_2] =\sum\limits_{p=2}^{m-2}
\gamma_{1,p}y_{p+1} - [y_2, x_2].$$

From the nilpotence condition of the superalgebra $L$ it follow
that in the decomposition of the product $[y_2, x_2]$ the basic
element $y_2$ does not participate, i.e. $[y_2, x_2] =
\sum\limits_{i=3}^{m} \gamma_{2,i}y_i.$ Therefore
$$[y_1, x_3] = ( \gamma_{1,2} - \gamma_{2,3}) y_3 + ( \gamma_{1,3}
- \gamma_{2,4}) y_4 + \dots + ( \gamma_{1,m-2} - \gamma_{2,m-1})
y_{m-1} - \gamma_{2,m}y_m,$$

In a similar way we obtain
$$[y_i, x_3] =( \gamma_{i, i+1} - \gamma_{i+1, i+2})y_{i+2} + \dots + (\gamma_{i, m-2} - \gamma_{i+1, m-1})
y_{m-1} - \gamma_{i+1, m}y_m, \ 2 \leq i \leq m-3.$$

Applying the above arguments for $[y_i, x_j], \ 4 \leq j \leq n$
we get that $[y_i, x_j]=\sum\limits_{j\geq i+2}(*)y_{j}$ for $4
\leq j \leq n$. Therefore without loss of generality one can
assume that the expression $\beta_{1,1}x_1 + \beta_{1,3}x_3 +
\dots + \beta_{1,n}x_n$ can be replaced by $x_1,$ i.e. we can
suppose $[y_1, y_1] =x_1.$ Consider the equalities
$$[x_1, y_1] = [[y_1, y_1], y_1] = \frac{1}{2}[y_1, [y_1, y_1]] = \frac{1}{2}[y_1, x_1]
=\frac{1}{2}y_2.$$

Since $[x_1, y_1] + [y_1, x_1] = \frac{3}{2}y_2 \in {\mathcal
R}(L)$ and $y_j = [y_j, x_1]$ for $1 \leq j \leq m_1-1,$ then
$y_3, y_4, \dots,y_{m_1} \in {\mathcal R}(L).$

Using induction and the following chain of equalities
$$[y_2, y_1] = [[y_1, x_1], y_1] = [y_1, [x_1, y_1]] + [[y_1, y_1], x_1] = \theta_1x_n,$$
$$[y_i, y_1] = [[y_{i-1}, x_1], y_1] = [y_{i-1}, [x_1, y_1]] + [[y_{i-1},
y_1], x_1] =[[y_{i-1}, y_1], x_1],$$
$$[x_j, y_1] = [[x_{j-1}, x_1], y_1] = [x_{j-1}, [x_1, y_1]] + [[x_{j-1}, y_1],
x_1]=[[x_{j-1}, y_1], x_1],$$ we establish that $$[y_i, y_1] =0 \
\mbox{for}\  3 \leq i \leq m_1$$ and $$[x_i, y_1] =
\alpha_{2,2}y_i+ \dots +\alpha_{2,m_1+2-i} y_{m_1} +
\alpha_{2,m_1+1} y_{m_1+i-1}+ \dots +$$ $$\alpha_{2,m_1 + m_2 + 2
- i} y_{m_1+m_2}+\alpha_{2,m_1 + m_2 + 1} y_{m_1+m_2+i-1} + \dots
+ \alpha_{2,m+2-i}y_m, \ 3\leq j \leq n.$$

The obtained products lead to $y_2 \notin L^4$ and that the basic
element $y_{m_1+m_2+ \dots +m_{t_0}+1}$ is generated by the
 products $[y_j, x_2], \ 2 \leq j \leq m_1.$

Since $y_j \in  {\mathcal R}(L)$ for $2 \leq j \leq  m_1$ and the
other basic elements $y_{m_1+1}, \dots, y_m$ are generated by
 products of the form $[[y_j,x_2], \dots],$ where $2 \leq j \leq m_1,$
then $\{y_2, y_3, \dots, y_{m}\} \in {\mathcal R}(L).$ Since
$[[y_j,x_2],y_1]=[y_j,[x_2,y_1]]+[[y_j,y_1],x_2]=0$ and

$$[[[y_j,x_2],\dots,x_2],y_1]=[[[y_j,x_2],\dots],[x_2,y_1]]+$$
$$[[[[y_j,x_2],\dots],y_1],x_2]=[[[[y_j,x_2],\dots],y_1],x_2]=\dots=0,$$
then we easily obtain that $$[y_{m_1+m_2+ \dots +m_{t}+1}, y_1] =
0, \ \mbox{for} \ 1\leq t \leq k-1.$$ Inductively we get
$$[y_j, y_1] = 0 \ \mbox{for}\  2 \leq j \leq m,$$ and that $x_3$
does not lie in $L^4.$ Thus, we obtain that
$$L^4 \subseteq \{x_4, \dots, x_n, y_3, \dots, y_m\}.$$ Hence $dim
L^4 < n+m -4,$ but this contradict the condition that the nilindex
is equal to $n+m.$ Thus, in the three classes of Theorem \ref{t2}
we obtain a contradiction with the assumption that the
superalgebra $L$ has nilindex equal to $n+m$ and therefore the
assertion of the theorem is proved. $\hfill \square$

From Theorem \ref{t3} we can assume that $x_1$ and $y_1$ are
generators of the superalgebra $L.$

\begin{thm}\label{t4} Let $L$ be a Leibniz superalgebra from $Leib_{n,m}$ with characteristic sequence
equal to $(n-1, 1 | m_1, m_2, \dots, m_k)$ and let $\{x_1, \
y_1\}$ be generators of $L.$ Then the superalgebra $L$ has
nilindex less than $n+m.$
\end{thm} \textit{Proof.} Let $L$ be a superalgebra satisfying the conditions of
the theorem. Then $$L^2 = \{x_2, x_3, \dots, x_n, y_2, y_3, \dots,
y_m\}.$$ Since  $y_{m_1+\dots +m_t+1}  \in L^2$ for any $t$ ($1
\leq  t \leq k-1$), we can conclude that $(\alpha_{1, m_1+\dots
+m_t+1} , \alpha_{2, m_1+\dots +m_t+1}, \dots, \alpha_{n,
m_1+\dots +m_t+1}) \neq (0, 0, \dots, 0)$ for any $t$ ($1 \leq t
\leq k-1).$ But this means that for any $t$ ($1 \leq t \leq k-1)$
the basic element $y_{m_1+\dots +m_t+1}$ is generated by the
products $[x_i, y_1], \ 1 \leq i \leq n.$

As in the proof of Theorem \ref{t3} denote by $y_{m_1+\dots
+m_{t_0}+1}$ the basic element which first is absent among the
elements
$$\{y_{m_1+1}, y_{m_1+m_2+1}, \dots, y_{m_1+m_2+ \dots+ m_{k-1}+1}
\}$$ in descending lower sequence. Then $(\alpha_{1, m_1+\dots
+m_{t_0}+1} , \alpha_{2, m_1+\dots +m_{t_0}+1}, \dots, \alpha_{n,
m_1+\dots +m_{t_0}+1}) \neq (0, 0, \dots, 0)$. Let $f$ be the
natural number such that $\alpha_{f, m_1 + \dots + m_{t_0}+1} \neq
0$ and $\alpha_{k, m_1 + \dots + m_{t_0}+1} = 0$ for $f \leq k
\leq n.$

We shall prove that $f=n.$ Let us suppose the opposite, i.e. $f <
n.$ Then for the powers of descending lower sequences we have the
following:
$$L^s = \{x_f, \dots, x_n, y_r, \dots, y_{m_1}, y_{m_1+1}, \dots,
y_{m_1+ \dots +m_{t_0}}, y_{m_1+ \dots + m_{t_0}+1}, \dots,
y_m\},$$
$$L^{s+1} = \{x_{f+1}, \dots, x_n, y_r, \dots, y_{m_1}, y_{m_1+1}, \dots,
y_{m_1+ \dots + m_{t_0}}, y_{m_1+ \dots + m_{t_0}+1}, \dots,
y_m\},$$
$$L^{s+2} = \{x_{f+1}, \dots, x_n, y_r, \dots, y_{m_1}, y_{m_1+1}, \dots,
y_{m_1+ \dots + m_{t_0}}, y_{m_1+ \dots + m_{t_0}+2}, \dots,
y_m\}.$$ From these we have that the elements $\{y_r, \dots,
y_{m_1}, y_{m_1+1}, \dots, y_{m_1+ \dots + m_{t_0}}, y_{m_1+ \dots
+ m_{t_0}+2}, \dots, \\ y_m\}$ are obtained from the products
$[x_{i}, y_1], \ f+1\leq i \leq n.$

The elements $\{y_r, \dots, y_{m_1}, y_{m_1+1}, \dots, y_{m_1+
\dots + m_{t_0}}, y_{m_1+ \dots + m_{t_0}+2}, \dots, y_m\}$ belong
to $L^{s+3}$ (because $\{x_{f+1}, \dots, x_m\} \in L^{s+2}$) and
hence $x_{f+1} \notin L^{s+3}.$ Therefore in the decomposition
$$[y_{m_1+ \dots + m_{t_0}+1}, y_1] = \beta_{m_1+ \dots +
m_{t_0}+1, f+1}x_{f+1} + \dots + \beta_{m_1+ \dots + m_{t_0}+1,
n}x_n$$ we have $\beta_{m_1+ \dots + m_{t_0}+1, f+1} \neq 0.$

Consider the equalities
$$[x_f , [y_1, y_1]] = 2[[x_f, y_1], y_1]
= 2[\alpha_{f,r}y_r + \alpha_{f,r+1}y_{r+1} + \dots +\alpha_{f,
m_1+ \dots + m_{t_0}+1} y_{m_1+ \dots + m_{t_0}+1} + \dots $$
$$+ \alpha_{f,m}y_m, y_1] =2 \alpha_{f, m_1+ \dots + m_{t_0}+1}
\beta_{m_1+ \dots + m_{t_0}+1, f+1} x_{f+1} + \sum\limits_{i\geq
f+2}(*)x_i.$$

On the other hand$$[x_f, [y_1, y_1]] = [x_f, \beta_{1,2}x_2 +
\beta_{1,3}x_3 + \dots + \beta_{1,n}x_n] = \sum\limits_{j\geq
f+2}(*)x_j.$$

Comparing the coefficients of the basic elements we obtain
$$\alpha_{f, m_1+ \dots + m_{t_0}+1} \beta_{m_1+ \dots + m_{t_0}+1,
f+1} =0,$$ which contradicts the conditions: $\alpha_{f, m_1+
\dots + m_{t_0}+1} \neq 0, \quad \beta_{m_1+ \dots + m_{t_0}+1,
f+1} \neq 0.$ Thus, we get a contradiction with the assumption
that $f < n.$ Now we shall study the case where $f=n,$ i.e.
$\alpha_{n, m_1+ \dots + m_{t_0}+1} \neq 0.$ In this case for some
natural number $p$ we have
$$L^p = \{y_{m_1+1}, \dots,
y_{m_1+ \dots + m_{t_0}}, y_{m_1+ \dots +m_{t_0}+1}, \dots,
y_m\}.$$ It is clear that if $k \geq3,$ then $dim L^{p+1}- dim L^p
\geq 2$ and we have a contradiction with the nilindex condition.
If $k=2$ then the vector space generated by the elements
$<y_{m_1+1}, \dots, y_{m_1+ m_2}>$ forms an ideal of the
superalgebra $L.$ The quotient superalgebra
$\overline{L}=L/<y_{m_1+1}, \dots, y_{m_1+ m_2}>$ is also two
generated and $C(\overline{L})=(n-1, 1 | m_1).$ Now applying Lemma
3.4 from \cite{AOKh} we get a contradiction, which completes the
proof of the theorem. $\hfill\square$

Let us investigate the case where both generators lie in odd part
of the superalgebra $L$. The following theorem clears up the
situation in this case.
\begin{thm} Let $L=L_0\oplus L_1$ be a superalgebra from $Leib_{n,m}$ with
characteristic sequence equal to $(n-1, 1 | m_1, m_2, \dots,
m_k),$ where $m_1\geq 2, \ n\geq 3$ and let both generators lie in
$L_1.$ Then the superalgebra $L$ has nilindex less than $n+m.$
\end{thm}
\textit{Proof.} Since both generators of the superalgebra $L$ lie
in $L_1,$ they are linear combinations of the elements $\{y_1,
y_{m_1+1}, \dots , y_{m_1+m_2+\dots + m_{k-1}+1} \}.$ Without loss
of generality we may assume that $y_1$ and $y_{m_1+1}$ are
generators.

Let $L^{2t} = \{x_i, x_{i+1}, \dots, x_n, y_j, \dots, y_m\}$ for
some natural number $t$ and let $z$ be an arbitrary element such
that $z \in L^{2t} \backslash L^{2t+1}.$ Then $z$ is generated by
the products of even an number of generators. Hence $z\in L_0$ and
$L^{2t+1} = \{x_{i+1}, \dots, x_n, y_j, \dots, y_m\}.$ In a
similar way, having $L^{2t+1} = \{x_{i+1}, \dots, x_n, y_j, \dots,
y_m\}$ we obtain $L^{2t+2} = \{x_{i+1}, \dots, x_n, y_{j+1},
\dots, y_m\}.$

From the above arguments we conclude that $n = m-1$ or $n = m-2.$
\\ Let us consider powers of $L:$ $$L^2 = \{x_1, x_2, \dots,
x_n, y_2, y_3, \dots, y_{m_1}, y_{m_1+2}, \dots, y_m\},$$
$$L^3 = \{A_1x_1 +A_2x_2, x_3, \dots, x_n, y_2, y_3, \dots, y_{m_1},
y_{m_1+2}, \dots, y_m\},$$
$$L^4 \supseteq \{A_1x_1 +A_2x_2, x_3, \dots, x_n, y_3, y_4, \dots, y_{m_1},
y_{m_1+3}, \dots, y_m\}.$$ Applying the above arguments we get
 that an element from the set $\{y_3, y_4, \dots, y_{m_1},
y_{m_1+3}, \\ \dots, y_m\}$ disappears in $L^4.$ If necessary then
by a shifting of basic elements we can suppose that $y_2 \notin
L^4$ without loss of the generality. Then
$$ L^4  = \{A_1x_1 +A_2x_2, x_3, \dots, x_n, y_3, y_4, \dots, y_{m_1},
y_{m_1+2}, \dots, y_m\},$$
$$ L^5  = \{x_3, \dots, x_n, y_3, y_4, \dots, y_{m_1},
y_{m_1+2}, \dots, y_m\}.$$

From these restrictions on the powers of $L$ in the following
products
$$[y_1, y_1 ] = \beta_{1,1}x_1 + \beta_{1,2}x_2 + \dots +
\beta_{1,n}x_n,$$
$$[y_2, y_1 ] = \beta_{2,2}(A_1x_1 +A_2x_2) + \beta_{2,3}x_3 + \dots + \beta_{2,n}x_n,$$
$$[y_1, y_{m_1+1}] = \gamma_{1,1}x_1 + \gamma_{1,2}x_2 + \dots + \gamma_{1,n}x_n,$$
$$[y_2, y_{m_1+1}] = \gamma_{2,2}(A_1x_1 +A_2x_2) + \gamma_{2,3}x_3 + \dots +
\gamma_{2,n}x_n,$$ we obtain the condition $(\beta_{2,2},
\gamma_{2,2}) \neq (0, 0).$

Let us introduce the notations
$$[x_1, y_1] = \alpha_{1,2}y_2 + \alpha_{1,3}y_3 + \dots + \alpha_{1,m_1} y_{m_1}
+ \alpha_{1,m_1+2} y_{m_1+2}  + \dots +  \alpha_{1,m} y_m,$$
$$[x_2, y_1] = \alpha_{2,2}y_2 + \alpha_{2,3}y_3 + \dots + \alpha_{2,m_1} y_{m_1}
+ \alpha_{2,m_1+2} y_{m_1+2}  + \dots +  \alpha_{2,m} y_m.$$
Consider the equalities
$$[x_1, [y_1, y_1]] =
2 [[x_1, y_1], y_1] = 2 [ \alpha_{1,2}y_2 + \alpha_{1,3}y_3 +
\dots + \alpha_{1,m_1} y_{m_1} + \alpha_{1,m_1+2} y_{m_1+2} +$$ $$
\dots +\alpha_{1,m} y_m, y_1] =2 \alpha_{1,2}\beta_{2,2}(A_1x_1
+A_2x_2) + \sum\limits_{i\geq 3}(*)x_i.$$ On the other hand
$$[x_1, [y_1, y_1]] = [x_1, \beta_{1,1}x_1 + \beta_{1,2}x_2 +
\dots + \beta_{1,n}x_n] = \sum\limits_{j\geq 3}(*)x_j.$$ Comparing
the coefficients of the basic elements in these equations we
obtain
$$\alpha_{1,2}\beta_{2,2} = 0.$$
Consider the product
$$[y_1, [y_1, x_1]] = [[y_1, y_1], x_1] - [[y_1, x_1], y_1] = [\beta_{1,1}x_1 +
\beta_{1,2}x_2 + \dots + \beta_{1,n}x_n, x_1] -$$ $$[y_2,y_1] =
-\beta_{2,2}(A_1x_1 +A_2x_2)+ \sum\limits_{s\geq 3}(*)x_s.$$ Since
$[y_1, [y_1, x_1]] = [y_1, y_2]$ then $[y_1, y_2] = -
\beta_{2,2}(A_1x_1 +A_2x_2) +\sum\limits_{s\geq 3}(*)x_s.$ \\From
the following the chain of equalities
$$[y_1, [x_1, y_1]]=[[y_1, x_1], y_1] - [[y_1, y_1], x_1] = [y_2, y_1] -
[\beta_{1,1}x_1 + \beta_{1,2}x_2 + $$ $$\dots +\beta_{1,n}x_n,
x_1]= \beta_{2,2}(A_1x_1 +A_2x_2) + \sum\limits_{t\geq 3}(*)x_t$$
and
$$[y_1, [x_1, y_1]] = [y_1,  \alpha_{1,2}y_2 + \alpha_{1,3}y_3 +
\dots + \alpha_{1,m_1} y_{m_1} + \alpha_{1,m_1+2} y_{m_1+2}  +
\dots +  \alpha_{1,m} y_m]=$$ $$ = - \alpha_{1,2}
\beta_{2,2}(A_1x_1 +A_2x_2) + \sum\limits_{p\geq 3}(*)x_p$$ we
obtain the restriction $\beta_{2,2} = -  \alpha_{1,2}
\beta_{2,2}.$ \\ Taking into account the condition
$\alpha_{1,2}\beta_{2,2} = 0$ we get $\beta_{2,2} = 0.$ Consider
$$ [y_1, [y_{m_1+1}, x_1]] = [[y_1,y_{m_1+1}], x_1] - [[y_1,
x_1],y_{m_1+1}]=[\gamma_{1,1}x_1 + \gamma_{1,2}x_2 + \dots
$$ $$ + \gamma_{1,n}x_n, x_1] -[y_2, y_{m_1+1}] =
-\gamma_{2,2}(A_1x_1 +A_2x_2) + \sum\limits_{q\geq 3}(*)x_q.$$
 On the other hand we have $ [y_1, [y_{m_1+1}, x_1]] = \sum\limits_{l\geq 3}(*)x_l.$
Comparing the coefficients of the basic elements we get
$\gamma_{2,2} = 0,$ which contradicts the condition $(\beta_{2,2},
\gamma_{2,2}) \neq (0, 0).$ Hence, we have $L^3  = \{x_3, \dots,
x_n, y_2, y_3, \dots, y_{m_1}, y_{m_1+2}, y_{m_1+3}, \dots,
y_m\},$ i.e. $A_1x_1+A_2x_2\notin L^3.$ Therefore the nilindex of
the superalgebra $L$ is less than $n+m.$ $\hfill\square$

The investigation of the cases where $L$ is a Leibniz superalgebra
with characteristic sequence $(n-1,1 | m_1, m_2, \dots, m_k),$
where either $m_1<2$ or $n<4,$ give us the same result.
Considering these cases consist is a simple routine work, mainly
repeating  the above technique, hence we omit details for these
cases.

\textbf{Acknowledgments.} \emph{The second author would like to gratefully
acknowledge the hospitality of the $\,$ "Institut f\"{u}r
Angewandte Mathematik",$\,$ Universit\"{a}t Bonn (Germany), as well as
the support by the DAAD foundation and by
the Nato-Reintegration Grant CBP.EAP.RIG.983169.}

\end{document}